# PERIODIC COPOLYMERS AT SELECTIVE INTERFACES: A LARGE DEVIATIONS APPROACH

By Erwin Bolthausen and Giambattista Giacomin

*Universität Zürich and Université Paris 7*

We analyze a $(1+1)$-dimension directed random walk model of a polymer dipped in a medium constituted by two immiscible solvents separated by a flat interface. The polymer chain is heterogeneous in the sense that a single monomer may energetically favor one or the other solvent. We focus on the case in which the polymer types are periodically distributed along the chain or, in other words, the polymer is constituted of identical stretches of fixed length. The phenomenon that one wants to analyze is the *localization at the interface*: energetically favored configurations place *most of* the monomers in the preferred solvent and this can be done only if the polymer sticks close to the interface.

We investigate, by means of large deviations, the energy–entropy competition that may lead, according to the value of the parameters (the strength of the coupling between monomers and solvents and an asymmetry parameter), to localization. We express the free energy of the system in terms of a variational formula that we can solve. We then use the result to analyze the phase diagram.

## 1. Introduction and results.

1.1. *The model.* Let $S = \{S_x\}_{x=0,1,\ldots}$ be a simple symmetric random walk starting at 0, that is, $S_0 = 0$, $S_x = \sum_{j=1}^{x} Y_j$, where $\{Y_j\}_{j=1,2,\ldots}$ are i.i.d. variables such that $\mathbb{P}(Y_1 = \pm 1) = 1/2$. For a *deterministic* sequence $\omega = \{\omega_x\}_{x \in \mathbb{N}}$, $\omega_x \in \{-1, 1\}$ and a parameter $h \geq 0$, we introduce the *Hamiltonian*

$$(1.1) \qquad H_{N,\omega,h}(S) \stackrel{\text{def}}{=} \sum_{x=1}^{N} (\omega_x + h) \operatorname{sign}(S_x),$$

and the probability measure $\mathbb{P}_{N,\omega,\lambda,h}$

$$\frac{d\mathbb{P}_{N,\omega,\lambda,h}}{d\mathbb{P}}(S) \stackrel{\text{def}}{=} \frac{\exp(\lambda H_{N,\omega,h}(S))}{Z_{N,\omega,\lambda,h}}, \tag{1.2}$$

where $\lambda \geq 0$ is the coupling constant, and $Z_{N,\omega,\lambda,h}$ is the normalization

$$Z_{N,\omega,\lambda,h} \stackrel{\text{def}}{=} \mathbb{E}[\exp(\lambda H_{N,\omega,h}(S))]. \tag{1.3}$$

We make the convention that if $S_x = 0$, then $\text{sign}(S_x) = \text{sign}(S_{x-1})$, which amounts to assigning the sign to the bond $((x-1, S_{x-1}), (x, S_x))$ rather than to the vertex $(x, S_x)$; see Figure 1. Our basic assumption on the sequence $\omega = \{\omega_x\}_{x \in \mathbb{N}}$ is that it is *centered* and *periodic*; that is, there exists $T \in \mathbb{N}$ such that $\omega_x = \omega_{x+2T}$ for every $x \in \mathbb{N}$ and $\sum_{x=1}^{2T} \omega_x = 0$. We write $T_\omega$ for the smallest such $T$. We exclude from our analysis the *trivial* $\omega$'s for which $\omega_{2k-1}\omega_{2k} = -1$ for every $k \in \mathbb{N}$ and in particular the case of $T_\omega = 1$, that is, $\omega_x = \pm(-1)^x$. In this case, $\sum_{x=1}^{N} \omega_x \text{sign}(S_x)$ is $0$ or $\pm 1$, due to the 2-periodicity of the random walk, and therefore, the influence of $\omega$ on the path measure is asymptotically negligible. We write $\mathcal{T}$ for the set of centered periodic sequences $\omega$ which are nontrivial. This is the *periodic version* of the random model considered in [3], where $\{\omega_x\}_x$ is a typical realization of a sequence of i.i.d. centered variables taking values $\pm 1$. We will often drop the dependence on $\lambda, h$ for notational convenience.

$\mathbb{P}_{N,\omega}$ is our model for a heterogeneous polymer near an interface, the $x$-axis, separating two media that interact with the monomers according to their $\omega$-characteristics. Possibly enlightening is the analogy with an oil/water interface and monomers, that is, $(S_{x-1}, S_x)$, that are either hydrophobic ($\omega_x = +1$) or hydrophilic ($\omega_x = -1$).

The free energy of such a model in the infinite volume limit is defined as

$$f_\omega(\lambda, h) = \lim_{N \to \infty} \frac{1}{N} \log Z_{N,\omega}. \tag{1.4}$$

The limit in (1.4) is easily seen to exist. We omit a direct proof as it follows from our more precise results, see Proposition 1.4. (An elementary direct proof can be given as in [3]; see, e.g., [7].)

1.2. *Localization and delocalization regions and the critical line.* An elementary estimate is

$$f_\omega(\lambda, h) \geq \lambda h, \tag{1.5}$$

for every $\omega$, $\lambda$ and $h$. This follows by observing that if $\Omega_N^+ \stackrel{\text{def}}{=} \{S : S_x \geq 0 \text{ for every } x = 1, 2, \ldots, N\}$, then

$$\frac{1}{N} \log Z_{N,\omega} \geq \frac{1}{N} \log \mathbb{E}\left[\exp\left(\lambda \sum_{x=1}^{N} (\omega_x + h) \text{sign}(S_x)\right); \Omega_N^+\right]$$

$$\text{(1.6)} \qquad = \frac{\lambda}{N} \sum_{x=1}^{N} (\omega_x + h) + \frac{1}{N} \log \mathbb{P}(\Omega_N^+)$$

$$= \lambda h + O\left(\frac{\log N}{N}\right), \qquad N \to \infty,$$

by the well-known estimate $1/\mathbb{P}(\Omega_N^+) = O(N^{1/2})$ (cf. [6], Chapter 3).

As in [3], motivated by the steps in (1.6), we partition the parameter space $\{(\lambda, h) : \lambda \geq 0, h \geq 0\}$ into two regions, $\mathcal{L}$ and $\mathcal{D}$:

(i) the *localized* region $\mathcal{L} = \{(\lambda, h) : f_\omega(\lambda, h) > \lambda h\}$;
(ii) the *delocalized* region $\mathcal{D} = \mathcal{L}^c = \{(\lambda, h) : f_\omega(\lambda, h) = \lambda h\}$.

We will discuss in Section 1.7 why and how (de)localization in the free energy sense is *equivalent* to pathwise de(localization). We set

$$\text{(1.7)} \qquad \phi_\omega(\lambda, h) = f_\omega(\lambda, h) - \lambda h.$$

The first result is:

PROPOSITION 1.1. *For every $\omega \in \mathcal{T}$ there exists a continuous nondecreasing function $h_c : [0, \infty) \longrightarrow [0, 1)$ such that $\mathcal{L} = \{(\lambda, h) : h < h_c(\lambda)\}$. Moreover, $h_c(0) = 0$ and $\lim_{\lambda \to \infty} h_c(\lambda) = 1$.*

We are going to focus on getting precise estimates on the *critical curve* $h_c(\cdot)$.

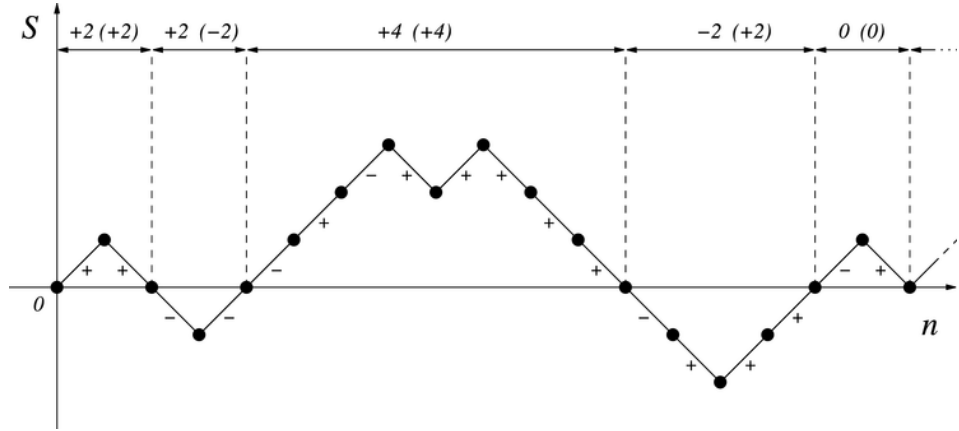

FIG. 1. *A polymer path: the quantities appearing on top of each excursion are $\sum_{x=j+1}^{k} \omega_x \operatorname{sign}(S_x)$ ($\sum_{x=j+1}^{k} \omega_x$), with $j$ and $k$, respectively, the beginning and the end point of the excursion.*

1.3. *A formula for the free energy.* We now give a series of definitions that lead to a rather explicit formula for the free energy. In Section 1.6 one can find an outline of the proof; this yields some intuition on this formula.

The basic idea of our approach is to split the Hamiltonian (1.1) as a sum of the contributions coming from the excursions of the random walk. An excursion is a portion $S_{2a} = 0, S_{2a+1}, \ldots, S_{2b-1}, S_{2b} = 0$ of the walk, where $S_x \neq 0$ for $2a < x < 2b$. Evidently, the relevant contribution coming from $\omega$ to one excursion is $\pm \sum_{x=2a+1}^{2b} \omega_x$, and this depends on $a$ and $b$ only through their values modulo $T_\omega$. It is therefore natural to define the following matrix indexed by the Abelian group $\mathbb{S} \stackrel{\text{def}}{=} \mathbb{Z}/T_\omega \mathbb{Z}$: for $\alpha, \beta \in \mathbb{S}$

$$\xi_{\alpha,\beta} \stackrel{\text{def}}{=} \sum_{x=2a+1}^{2b} \omega_x, \tag{1.8}$$

which is well defined by choosing representatives $a \in \alpha$ and $b \in \beta$ with $a < b$. Evidently

$$\xi_{\alpha,\beta} = \xi_{0,\beta} - \xi_{0,\alpha}. \tag{1.9}$$

We set $\xi_\star \stackrel{\text{def}}{=} \max_{\alpha,\beta} |\xi_{\alpha,\beta}| \leq 2T_\omega$. Notice that for every $\omega \in \mathcal{T}$ the matrix $\{\xi_{\alpha,\beta}\}_{\alpha,\beta}$ is not identically zero.

As an example, take the simplest case of an $\omega$-sequence: $+, +, -, -, +, +, -, -, \cdots$ ($T_\omega = 2$), which yields $\xi = \begin{pmatrix} 0 & 2 \\ -2 & 0 \end{pmatrix}$.

We need also some notation for the random walk $S$. Let $\eta$ be the first return time to 0, that is, $\eta \stackrel{\text{def}}{=} \inf\{x \geq 1 : S_x = 0\}$, and set $K(x) \stackrel{\text{def}}{=} \mathbb{P}(\eta = x)$ for $x \in 2\mathbb{N}$, $p_{\alpha,\beta} \stackrel{\text{def}}{=} \mathbb{P}(\eta/2 \in \beta - \alpha)$ for $\alpha, \beta \in \mathbb{S}$, and $K_{\alpha,\beta}(x) \stackrel{\text{def}}{=} \mathbb{P}(\eta = x | \eta/2 \in \beta - \alpha)$. Note that $\{p_{\alpha,\beta}\}_{\alpha,\beta}$ is bistochastic. As it is well known,

$$\lim_{x \in 2\mathbb{N}, x \to \infty} x^{3/2} K(x) = \sqrt{2/\pi} =: \mathsf{C}_K, \tag{1.10}$$

see, for example, [6], Chapter 3.

We set for $x \in \mathbb{N}$ and $\alpha, \beta \in \mathbb{S}$

$$\Phi_{\alpha,\beta}^{\lambda,h}(x) \stackrel{\text{def}}{=} \log\left(\frac{1 + \exp(-2(\lambda \xi_{\alpha,\beta} + \lambda h x))}{2}\right), \tag{1.11}$$

and in turn the $\mathbb{S} \times \mathbb{S}$ matrix with positive entries

$$A_{\alpha,\beta} = A_{\alpha,\beta}(b, \lambda, h) \stackrel{\text{def}}{=} p_{\alpha,\beta} \sum_x K_{\alpha,\beta}(x) \exp(\Phi_{\alpha,\beta}^{\lambda,h}(x) - bx), \tag{1.12}$$

with $b \geq 0$. Let us denote by $Z = Z(b, \lambda, h)\ (>0)$ the Perron–Frobenius (maximal) eigenvalue. We observe that $A_{\alpha,\beta}(\cdot, \lambda, h)$ is a decreasing function for every $\alpha$ and $\beta$. This implies that $Z(\cdot, \lambda, h)$ is decreasing (cf. [10], Chapter 1). For the same reason $Z(b, \lambda, \cdot)$ is decreasing, too. We have the following:

THEOREM 1.2. *Denote by $\tilde{b} = \tilde{b}(\lambda, h)$ the unique solution of $Z(\tilde{b}, \lambda, h) = 1$, if such a solution exists, and set $\tilde{b}(\lambda, h) = 0$ if such a solution does not exist. We have that*

$$\phi_\omega(\lambda, h) = \tilde{b}(\lambda, h), \tag{1.13}$$

*for every $\lambda$ and $h$.*

An immediate consequence of this formula, of the fact that $Z(b, \lambda, \cdot)$ is decreasing and of Proposition 1.1 is that $h_c(\lambda)$ is uniquely determined by the equation $Z(0, \lambda, h_c(\lambda)) = 1$.

1.4. *Estimates on the critical line.* From Theorem 1.2 we extract the precise asymptotic behavior of the critical curve at small values of $\lambda$.

THEOREM 1.3. *For every $\omega \in \mathcal{T}$, as $\lambda \searrow 0$ we have that*

$$h_c(\lambda) = m_\omega \lambda^3 (1 + o(1)), \tag{1.14}$$

*where*

$$m_\omega = \left( \frac{1}{2T_\omega} \sum_{\alpha, \beta} p_{\alpha, \beta} \xi_{\alpha, \beta}^2 \right)^2. \tag{1.15}$$

*Moreover, there exists a positive constant $M_\omega$ such that as $\lambda \nearrow \infty$,*

$$h_c(\lambda) = 1 - (M_\omega + o(1))\lambda^{-1}. \tag{1.16}$$

In Section 5.2 one can find an expression for $M_\omega$.

1.5. *On related copolymer models.* A large amount of papers dealing with periodic copolymers can be found in the literature (mostly in the area of chemistry and physics). We single out some of them, divided into two categories: results about systems with period 2 ($T_\omega = 1$) and results about more general copolymers.

*Period-2 copolymers.* We stress once again that the $T_\omega = 1$ case leads to a trivial model in our case. This is due to an evident cancellation connected to the fact that the $\omega$-periodicity coincides with the periodicity of the walk. One may, however, modify *slightly* the definition of the model by taking a different convention for the sign of zero and the situation may change, leading to a localization–delocalization phenomenon. We mention in particular the work [12] in which $\omega_x = (-1)^{x+1}$ and $\text{sign}(0) = -1$. Since this can be interpreted as choosing the interface line at height $1/2$, we will refer to this model as the $1/2$-interface (copolymer) model. The authors point out in particular that such a model, under *one-step decimation*, becomes a homogeneous polymer model at a *penetrable attractive interface*, which is known to be *exactly solvable* (we refer also to [9] and references therein for exact computations on this model). In probabilistic language, *one-step decimation* means simply to consider the marginal of the polymer measure over odd (or even) sites. The computation is elementary and the arising model is simply a random walk that prefers the upper half plane if $h > 0$ and that receives a reward each time that crosses the interface. We mention also the result [15] in which $\text{sign}(0) = 0$, but in which the asymmetry is replaced by a penalization for the walk to touch zero; in this case, again by one-step decimation, the system reduces to the walk with a reward (positive or negative) at the origin. We signal also the complete analysis obtained in [8] for a Gaussian random walk with alternating $\omega$.

Some of the results in [12] may at first look to be in contradiction with ours: in particular, that in [12] it is shown that $h_c(\lambda) = (1 + o(1))\lambda$ as $\lambda \searrow 0$. We point out, however, that if one introduces the $h$ asymmetry for the model with $\text{sign}(0) = 0$ (cf. [15]), one can show that $h_c(\lambda) = (1 + o(1))\lambda^3$ (we refer to this model as the *neutral interface* model). It is, however, not too difficult to understand the mechanism that leads to the different phenomenology in the 1/2-interface model:

(a) in the 1/2-interface model *any excursion is favorable*, in the sense that, if one considers successive crossings of the 1/2-interface, each time the energetic gain from $\omega$ amounts to $+1$;

(b) in the neutral interface model this does not happen: crossing the interface is not enough to get a positive energy contribution or, in other words, there are favorable and unfavorable excursions.

The two results that we have just mentioned can be established also via our approach; in the $T_\omega = 1$ case the arising variational problem can be written in a fully explicit fashion. This of course leads to results approaching the completeness of the analysis in [12].

*More general periodic polymers.* Among the physics papers on the case $T_\omega > 1$ we single out the one of Sommer and Daoud [16], who consider the case in which $\omega$ is made of alternating blocks of length $T_\omega$ of the same sign (this model is referred to as *diblocks* model). While no mathematically precise model is given, the authors argue, on the base of scaling arguments and of a renormalization group analysis, in favor of $h_c(\lambda) \sim T_\omega^3 \lambda^3$. This agrees with our result not only because of the correct $\lambda$ dependence, but also because $m_\omega$ behaves like $T_\omega^3$ for large $T_\omega$; it must be noted, however, that this behavior is restricted to the diblocks case and it is rather easy to see that:

(a) fastest growth of $m_\omega$ in $T_\omega$ is obtained for the diblocks model;

(b) one can construct any intermediate behavior down to the extreme case of $m_\omega = O(1/T_\omega)$ as $T_\omega$ grows. For example, the latter case is achieved by starting with the periodic (forbidden!) configuration $+-+-+-\cdots$, and switching $-+ \longrightarrow +-$ at regular intervals (say $2k$ sites), obtaining thus an element $\omega \in \mathcal{T}$ with $T_\omega = k$ and satisfying the desired property.

We signal also that precise large-period developments for the related problem of copolymer *adsorption*, for a very special family of $\omega$'s, have been established by combinatorial methods in [13].

Moreover, it is interesting to recall that if $\omega$ is a typical realization of an i.i.d. sequence of centered random variables, the phase diagram is still split in two regions $\mathcal{D}$ and $\mathcal{L}$ by a continuous function $h_c(\cdot)$ for which a strict analog of Proposition 1.1 holds (cf. [3]). However, in that case the derivative of $h_c(\cdot)$ in zero (exists and it) is positive and bounded above by 1 (cf. [3]; see also [11] for various physical predictions). The appearance of a positive slope in the random $\omega$ case is nontrivial: an excursion of length $L$, $L$ large, leads to an energetic gain of the order of $\sqrt{L}$, if the sign of the excursion agrees with the sign of the $\omega$ fluctuation. This effect is of course not present in the periodic case, in which the energetic gain is always $O(1)$. As we have seen, a positive slope phenomenon may (and will) be observed if the interface is not neutral and attractive.

1.6. *Sketch of the proofs: a sharp energy–entropy argument.* As already pointed out, our arguments start with expressing the energy of the copolymer in terms of the return times to zero of $S$. These are defined recursively by

$$(1.17) \qquad \eta_0 \stackrel{\text{def}}{=} 0 \quad \text{and} \quad \eta_{k+1} \stackrel{\text{def}}{=} \inf\{x > \eta_k : S_x = 0\} \qquad \text{for } k \geq 1.$$

Of course once the sequence $\{\eta_k\}_k$ is fixed, to determine the energy we need to know the sign of $S_x$ for $x \in \{\eta_k + 1, \ldots, \eta_{k+1}\}$. However, the sequence of signs is i.i.d. and centered, so that these degrees of freedom are easily integrated out; see the beginning of Section 3 for the straightforward details.

It turns out to be practical to consider the sequence $\{[\eta_k/2], [\eta_{k+1}/2], \eta_{k+1} - \eta_k\}_{k=1,2,\ldots}$, where, for $x \in \mathbb{N}$, $[x] = [x]_\mathbb{S}$ denotes the class in $\mathbb{S}$. This sequence is a Markov chain in $\mathbb{S} \times \mathbb{S} \times 2\mathbb{N}$ with transition probabilities

$$(1.18) \qquad P((\alpha, \beta, x), (\alpha', \beta', x')) \stackrel{\text{def}}{=} \delta_{\beta, \alpha'} p_{\alpha', \beta'} K_{\alpha', \beta'}(x').$$

It is immediate to see that the stationary distribution $\pi_{\text{eq}}$ of this Markov chain is given by

$$(1.19) \qquad \pi_{\text{eq}}(\alpha, \beta, x) \stackrel{\text{def}}{=} \frac{1}{T_\omega} p_{\alpha, \beta} K_{\alpha, \beta}(x).$$

Since the energy can be essentially expressed in terms of the empirical measure of this Markov chain [cf. (3.4)], that is, the frequency with which each value of $(\alpha, \beta, x)$ is observed in the Markov chain sequence, one can in turn express the leading asymptotics of the exponential of the energy via a variational formula. This formula evaluates the competition between the energy and the entropy of the system, and the latter is the large deviation functional for the empirical measure.

More in detail, if $\mu$ is a probability measure on $\mathbb{S} \times \mathbb{S} \times 2\mathbb{N}$, we denote by $\mu_1, \mu_2, \mu_3$ the marginals on $\mathbb{S}$ and $2\mathbb{N}$, respectively. Let $\mathcal{P}$ be the set of probability measures $\mu$ on $\mathbb{S} \times \mathbb{S} \times 2\mathbb{N}$ satisfying $\mu_1 = \mu_2$, and $\mu(\{(\alpha, \beta, x) : x/2 \in \beta - \alpha\}) = 1$. Note that $\pi_{\text{eq}} \in \mathcal{P}$.

We set now

$$(1.20) \quad I(\mu) \stackrel{\text{def}}{=} \begin{cases} \sum_{\alpha, \beta, x} \mu(\alpha, \beta, x) \log\left(\frac{\mu(\alpha, \beta, x)}{\mu_1(\alpha) p_{\alpha, \beta} K_{\alpha, \beta}(x)}\right), & \text{if } \mu \in \mathcal{P}, \\ +\infty, & \text{if } \mu \notin \mathcal{P}, \end{cases}$$

with $0 \log 0 = 0$. If we put

$$(1.21) \quad Q(\mu) \stackrel{\text{def}}{=} \sum_{\alpha, \beta, x} \Phi_{\alpha, \beta}^{\lambda, h}(x) \mu(\alpha, \beta, x) - I(\mu),$$

then we have the following variational formula for $\phi_\omega$:

PROPOSITION 1.4. *For every $\omega \in \mathcal{T}$, $\lambda \geq 0$ and $h \geq 0$ we have*

$$(1.22) \quad \phi_\omega(\lambda, h) = \sup_{t > 0} t \sup\left\{Q(\mu) : \mu \in \mathcal{P}, \sum_x x \mu_3(x) \leq 1/t\right\}.$$

REMARK 1.5. We observe that, by (1.9), $\sum_{\alpha, \beta, x} \xi_{\alpha, \beta} \mu(\alpha, \beta, x) = 0$ for $\mu \in \mathcal{P}$ and it will turn out to be more practical at times to change the definition of $\Phi_{\alpha, \beta}^{\lambda, h}(x)$ by adding $\lambda \xi_{\alpha, \beta}$ so that

$$(1.23) \quad \Phi_{\alpha, \beta}^{\lambda, h}(x) = \psi(\lambda \xi_{\alpha, \beta} + \lambda h x) - \lambda h x,$$

with $\psi(\cdot) = \log \cosh(\cdot)$, and Proposition 1.4 still holds.

The variational problem in Proposition 1.4 can be solved (almost) explicitly. In order to explain how the solution looks, let us construct, in a fairly standard way, a perturbation of the transition probabilities (1.18) by defining the following family of functions on $\mathbb{S} \times \mathbb{S} \times 2\mathbb{N}$:

$$(1.24) \quad \mathcal{A}(\alpha, \beta, x) = \mathcal{A}_{b, \lambda, h}(\alpha, \beta, x) \stackrel{\text{def}}{=} p_{\alpha, \beta} K_{\alpha, \beta}(x) \exp(\Phi_{\alpha, \beta}^{\lambda, h}(x) - bx),$$

where $b \geq 0$ is a parameter. Of course $\sum_x \mathcal{A}_{b, \lambda, h}(\alpha, \beta, x)$ coincides with $A_{\alpha, \beta}(b, \lambda, h)$. Let us denote by $\{v_\alpha\}_{\alpha \in \mathbb{S}}$ the unique (up to scaling) positive right eigenvector of $A$ with eigenvalue $Z$ and by $\{\pi(\alpha)\}_{\alpha \in \mathbb{S}}$ the normalized left eigenvector of $\{A_{\alpha,\beta} v_\beta / v_\alpha\}_{\alpha,\beta}$. Then the measure $\mu_b^{\lambda,h}$ on $\mathbb{S} \times \mathbb{S} \times 2\mathbb{N}$ defined by

$$(1.25) \qquad \mu_b^{\lambda,h}(\alpha, \beta, x) \stackrel{\text{def}}{=} \frac{1}{Z} \pi(\alpha) \mathcal{A}(\alpha, \beta, x) \frac{v_\beta}{v_\alpha}$$

is in $\mathcal{P}$.

We are going to show (see Lemma 2.2 in Section 2) that the supremum over $\mu \in \mathcal{P}$ in (1.22) is attained in the set $\{\mu_b^{\lambda,h} : b > 0\}$. The computation of $\phi_\omega(\lambda, h)$ boils down, therefore, to optimizing over $b$. This can, once again, be done explictly, obtaining thus Theorem 1.2.

1.7. *A (quick) look at pathwise results.* For completeness we point out here that, if $(\lambda, h) \in \mathcal{L}$, for every $\varepsilon > 0$ there exists $C > 0$ such that

$$(1.26) \qquad \mathbb{P}_{N,\omega,\lambda,h}(|S_x| > L) \leq C \exp(-(\phi_\omega(\lambda, h) - \varepsilon)L),$$

for every $N$, every $x \leq N$ and every $L > 0$. This result can be easily extracted by applying the technique in [14]; see [1] for further results on the localized phase of the random $\omega$ model. On the other side, since we are able to solve explicitly the variational problem associated to the free energy, if $(\lambda, h) \in \mathcal{L}$ by large deviation arguments we have detailed information on the empirical measure of the copolymer, which converges in probability, in the $N \to \infty$ limit, to $\mu_{\tilde{b}(\lambda,h)}^{\lambda,h}$. This gives in particular that if we set $\ell_N = \max\{k \in \mathbb{N} : \eta_k \leq N\}$, then $\ell_N / N$ tends to $1/\sum_{\alpha,\beta,x} x \mu_{\tilde{b}(\lambda,h)}^{\lambda,h}(\alpha, \beta, x)$, in $\mathbb{P}_{N,\omega,\lambda,h}(dS)$-probability, or, in other words,

$$(1.27) \qquad \lim_{N \to \infty} \frac{1}{\ell_N} \sum_{k=1}^{\ell_N} (\eta_k - \eta_{k-1}) = \sum_{\alpha,\beta,x} x \mu_{\tilde{b}(\lambda,h)}^{\lambda,h}(\alpha, \beta, x).$$

In particular, one sees that there is a continuous blow-up of the typical excursion length approaching the delocalization region.

On the other side, if $(\lambda, h)$ belongs to the interior of $\mathcal{D}$, one can show (see the last section of [1]) that for every $L > 0$

$$(1.28) \qquad \lim_{N \to \infty} \frac{1}{N} \sum_{x=1}^{N} \mathbb{P}_{N,\omega,\lambda,h}(S_x > L) = 1.$$

It should be possible to improve (1.28) strongly, leading in particular to the Brownian scaling results in [9].

The paper is organized as follows: in Section 2 we study the solutions of the variational problem in Proposition 1.4. We show in particular that the right-hand side of (1.13) coincides with the right-hand side of (1.22). In Section 3 we prove our basic variational formula, that is, Proposition 1.4, completing thus the proof of Theorem 1.2. In Section 4 we study the existence of a critical line for the model and we partly prove Proposition 1.1. The proof is completed in the last section, where Theorem 1.3 is established.

## 2. The variational problem.

2.1. *The linear algebra setup.* Throughout this work we will make repeated use of the results of the Perron–Frobenius theory (see, e.g., [5], Chapter 3). Let us quickly recall these facts. For $T \in \mathbb{N}$ let us denote by $\mathbb{M}^+(T)$ the set of $T \times T$ matrices with nonnegative entries which are irreducible. If $A \in \mathbb{M}^+(T)$, then there exists a unique right eigenvector $v = v(A) \in (0, \infty)^T$, normalized by $\sum_{i=1}^T v_i = \|v\|_1 = 1$. The corresponding left eigenvector is denoted by $\tilde{v}$, $\|\tilde{v}\|_1 = 1$. If we call $\rho_{\max} = \rho_{\max}(A)$ the associated (positive) eigenvalue (the *maximal eigenvalue*), then for any other root $\rho \in \mathbb{C}$ of the characteristic polynomial we have that $|\rho| \leq \rho_{\max}$ and $\rho_{\max}$ is a simple root.

If we have a function $A : I \longrightarrow \mathbb{M}^+(T)$, $I$ open interval of $\mathbb{R}$, then the regularity of $A(\cdot)$ passes directly to the Perron–Frobenius eigenvalue and eigenvectors; see, for example, [4], Chapters 2 and 4. In particular, if $A(\cdot)$ is differentiable, then $\rho_{\max}(A(\cdot))$ is differentiable, too, and

$$\frac{d}{dt}\rho_{\max}(A(t)) = \frac{\sum_{\alpha,\beta} A'_{\alpha,\beta}(t)\tilde{v}_\alpha(A(t))v_\beta(A(t))}{\sum_\alpha \tilde{v}_\alpha(A(t))v_\alpha(A(t))}. \tag{2.1}$$

2.2. *Solutions to the variational problem.* In order to lighten the exposition we modify (and, hopefully, simplify) somewhat the notation with respect to the Introduction. The major change is that for the marginals $\mu_1$, $\mu_2$ and $\mu_3$ of $\mu \in \mathcal{P}$ we are going to drop the numerical index. This abuse of notation is of course partly justified for the first two marginals, but in order to avoid misunderstandings the argument(s) of $\mu$ will always be explicitly given: so $\mu(\alpha)$ is $\mu_1(\alpha)$ and $\mu(x) = \mu_3(x)$. In the same way, $\mu(\alpha, \beta)$ is of course a notation for $\sum_x \mu(\alpha, \beta, x)$. With this convention [cf. with (1.25)]

$$\mu_b(\alpha, \beta, x) = \mu_b(\alpha) p_{\alpha,\beta} K_{\alpha,\beta}(x) \left( \frac{\exp(\Phi_{\alpha,\beta}(x) - bx)}{Z} \right) \left( \frac{v_\beta}{v_\alpha} \right). \tag{2.2}$$

Notice that we have dropped the dependence on $\lambda$ and $h$ in $\Phi_{\alpha,\beta}^{\lambda,h}$.

If at times we will drop the explicit dependence on one or more parameters, in other situations the opposite tendency will prevail and we will write $A(b, \lambda, h)$, $A(b, \lambda)$, $v(b, \lambda)$, $\mu_b^\lambda$, and so on.

Moreover, for the results of this section the details of $\Phi_{\alpha,\beta}^{\lambda,h}(x)$ are inessential: all we are going to use is the regularity of $\Phi$ with respect to $\lambda$ and $h$ and that $\sup_x |\Phi_{\alpha,\beta}^{\lambda,h}(x)| < \infty$ for every $\alpha, \beta, \lambda$ and $h$.

We have the following results:

LEMMA 2.1. *For every $b > 0$*

$$\frac{d}{db} \log Z(b) = -\sum_x x \mu_b(x). \tag{2.3}$$

PROOF. We start by observing that one can construct the Markov chain on $\mathbb{S} \times \mathbb{S} \times 2\mathbb{N}$ that makes transition from $(\alpha, \beta, x)$ to $(\alpha', \beta', x')$ with probability $\mu_b(\alpha', \beta', x')\delta_{\beta,\alpha'}/\mu_b(\alpha')$. Notice that $\mu_b$ is invariant for such a chain and we denote by $\{\alpha_j, \alpha_{j+1}, \Delta\eta_j\}_j$ the stationary process and by $\langle \cdot \rangle$ its expectation. Observe, moreover, that for every $k \in \mathbb{N}$,

$$1 = \sum_{\substack{\alpha_0,\ldots,\alpha_k \\ x_0,\ldots,x_{k-1}}} \mu_b(\alpha_0, \alpha_1, x_0) \frac{\mu_b(\alpha_1, \alpha_2, x_1)}{\mu_b(\alpha_1)} \cdots \frac{\mu_b(\alpha_{k-1}, \alpha_k, x_{k-1})}{\mu_b(\alpha_{k-1})}$$

$$(2.4) \quad = \frac{1}{Z(b)^k} \sum_{\substack{\alpha_0,\ldots,\alpha_k \\ x_0,\ldots,x_{k-1}}} \left[\frac{\mu_b(\alpha_0) v_{\alpha_k}}{v_{\alpha_0}}\right] \exp\left(\sum_{j=0}^{k-1} \Phi_{\alpha_j,\alpha_{j+1}}(x_j) - b\sum_{j=0}^{k-1} x_j\right)$$

$$\times \prod_{j=0}^{k-1} p_{\alpha_j,\alpha_{j+1}} K_{\alpha_j,\alpha_{j+1}}(x_j).$$

This formula gives a family of expressions for $Z(b)$, indexed by $k$. Notice that the term between brackets $[\cdots]$ depends on $b$, but it is bounded and its derivative with respect to $b$ is bounded, too. By taking the derivative with respect to $b$ and letting $k$ go to infinity one obtains

$$(2.5) \quad \frac{d}{db}\log Z(b) = -\lim_{k\to\infty} \frac{1}{k}\left\langle \sum_{j=0}^{k-1} \Delta\eta_j \right\rangle = -\sum_x x\mu_b(x).$$

The proof is therefore complete. An alternative proof of (2.3) can of course be extracted directly from (2.1). □

LEMMA 2.2. *Set $f(b) = \sum_x x\mu_b(x)$.*

(i) *$f \in C(0, \infty)$ is decreasing. Moreover, $\lim_{b \nearrow \infty} f(b) = 2$ and $\lim_{b \searrow 0} f(b) = \infty$.*

(ii) *For every $\mu \in \mathcal{P}$ we have that if $\sum_x x\mu(x) = \sum_x x\mu_b(x)$, then*

$$(2.6) \quad Q(\mu) \le Q(\mu_b),$$

*and equality holds only if $\mu = \mu_b$.*

PROOF. The proof is a bit indirect: we start by establishing the continuity and the limits claimed in (i). We then prove (ii), that will imply the (strict) monotonicity of $f$.

First of all we observe that $A_{\alpha,\beta}(\cdot) \in C^0(0, \infty)$ and therefore $Z(\cdot), v_\beta(\cdot)/v_\alpha(\cdot)$ and $\mu_\cdot(\alpha)$ are in $C^0(0, \infty)$. Recalling that $\Phi_{\alpha,\beta}$ is bounded, one obtains that $f(\cdot) \in C^0(0, \infty)$.

In order to deal with the limits of $f$ at the boundary of $(0, \infty)$, we first remark that the expression in the right-hand side of (1.12) makes sense for $b = 0$, thus defining $A(0) \in \mathbb{M}^+(T_\omega)$ and $\lim_{b \searrow 0} A(b) = A(0)$. This guarantees that $Z(b)$, $v_\alpha(b)$ and $\mu_b(\alpha)$ tend to finite limits as $b$ vanishes: since $\sum_x x K_{\alpha,\beta}(x) = \infty$ we conclude that $\lim_{b \searrow 0} f(b) = +\infty$. On the other side, it is rather immediate to see that $\lim_{b \nearrow \infty} \mu_b(\alpha, \beta, x) = 1/T_\omega$ if $\beta - \alpha = 1$ and $x = 2$, and $\lim_{b \nearrow \infty} \mu_b(\alpha, \beta, x) = 0$ otherwise. This yields $\lim_{b \nearrow \infty} f(b) = 2$.

For what concerns (2.6) we first note that

$$
\begin{aligned}
I(\mu) &= \sum_{\alpha,\beta,x} \left[ \log\left( \frac{\mu_b(\alpha,\beta,x)}{\mu_b(\alpha) p_{\alpha,\beta} K_{\alpha,\beta}(x)} \right) \right. \\
&\quad \left. + \log\left( \frac{\mu(\alpha,\beta,x)}{\mu_b(\alpha,\beta,x)} \right) + \log\left( \frac{\mu_b(\alpha)}{\mu(\alpha)} \right) \right] \mu(\alpha,\beta,x), \\
(2.7) \quad &= \sum_{\alpha,\beta,x} \Phi_{\alpha,\beta}(x) \mu(\alpha,\beta,x) - b \sum_x x \mu_b(x) - \log Z(b) \\
&\quad + \sum_{\alpha,\beta,x} \log\left( \frac{v_\beta(b)}{v_\alpha(b)} \right) \mu(\alpha,\beta,x) + \widetilde{H}(\mu|\mu_b),
\end{aligned}
$$

where $\widetilde{H}(\mu|\mu_b)$ is the difference of relative entropies

$$
(2.8) \quad \sum_{\alpha,\beta,x} \log\left( \frac{\mu(\alpha,\beta,x)}{\mu_b(\alpha,\beta,x)} \right) \mu(\alpha,\beta,x) - \sum_\alpha \log\left( \frac{\mu(\alpha)}{\mu_b(\alpha)} \right) \mu(\alpha).
$$

We claim that $\widetilde{H}(\mu|\mu_b) > 0$, unless $\mu = \mu_b$ (see below for a proof of this claim). Moreover, we observe that

$$
(2.9) \quad \sum_{\alpha,\beta,x} \log\left( \frac{v_\beta(b)}{v_\alpha(b)} \right) \mu(\alpha,\beta,x) = \sum_{\alpha,\beta} \log\left( \frac{v_\beta(b)}{v_\alpha(b)} \right) \mu(\alpha,\beta) = 0,
$$

since the marginals of $\mu(\alpha,\beta)$ are identical.

Applying these observations to (2.7), we obtain that for every $b$

$$
(2.10) \quad Q(\mu) = b \sum_x x \mu(x) + \log Z(b) - \widetilde{H}(\mu|\mu_b),
$$

and since by a straightforward computation

$$
(2.11) \quad Q(\mu_b^{\lambda,h}) = b \sum_x x \mu_b^{\lambda,h}(x) + \log Z(b, \lambda, h),
$$

Part (ii) of the statement is proven. Of part (i) of the statement we are left with establishing the strict monotonicity of $f$. Since $f$ is continuous, if $f$ is not strictly monotonic, there exists $b_1 < b_2$ such that $f(b_1) = f(b_2)$. By applying (2.10) with $b = b_2$ and (2.11) we obtain

$$
(2.12) \quad Q(\mu_{b_1}) = b_2 f(b_2) + \log Z(b_2) - \widetilde{H}(\mu_{b_1}|\mu_{b_2}) \leq Q(\mu_{b_2}).
$$

Of course, by exchanging $b_1$ and $b_2$ we obtain the reversed inequality, that is, $Q(\mu_{b_2}) \leq Q(\mu_{b_1})$, and therefore that $\widetilde{H}(\mu_{b_1}|\mu_{b_2}) = 0$, which, by the claim, forces $\mu_{b_1} = \mu_{b_2}$, which is clearly impossible. Therefore $f$ is strictly monotonic.

Let us now establish the claim that $\widetilde{H}(\mu|\mu_b)$ is positive unless $\mu = \mu_b$. From (2.8) we observe that $\widetilde{H}(\mu|\mu_b)$ can be viewed as an average of relative entropies: set $\mu(\beta, x|\alpha) = \mu(\alpha, \beta, x)/\mu(\alpha)$ whenever $\mu(\alpha) \neq 0$, so that

$$\begin{aligned}
\widetilde{H}(\mu|\mu_b) &= \sum_{\alpha,\beta,x} \log\left(\frac{\mu(\beta,x|\alpha)}{\mu_b(\beta,x|\alpha)}\right)\mu(\alpha,\beta,x) \\
&= \sum_\alpha \mu(\alpha) H(\mu(\cdot,\cdot|\alpha)|\mu_b(\cdot,\cdot|\alpha)),
\end{aligned}$$
(2.13)

where $H$ is the standard relative entropy. Notice that $\mu \ll \mu_b$ by definition of $\mathcal{P}$ and $\mu_b$. Therefore $\widetilde{H}$ is nonnegative. If $\widetilde{H}(\mu|\mu_b) = 0$, then $H(\mu(\cdot,\cdot|\alpha)|\mu_b(\cdot,\cdot|\alpha)) = 0$, and therefore $\mu(\cdot,\cdot|\alpha) = \mu_b(\cdot,\cdot|\alpha)$, for every $\alpha$ such that $\mu(\alpha) \neq 0$. This implies that $\mu(\alpha,\beta)/\mu(\alpha) = \mu_b(\alpha,\beta)/\mu_b(\alpha)$ for $\mu(\alpha) \neq 0$ and therefore $\sum_\alpha \mu(\alpha)(\mu_b(\alpha,\beta)/\mu_b(\alpha)) = \mu(\beta)$. If we set $A_{\alpha,\beta} = \mu_b(\alpha,\beta)/\mu_b(\alpha)$, we see that $A$ is a stochastic matrix in $\mathbb{M}^+(T_\omega)$; $\mu$ is therefore the unique nonnegative normalized left eigenvector. It is, however, immediate to verify that $\mu_b$ is also a nonnegative normalized left eigenvector; therefore $\mu(\alpha) = \mu_b(\alpha) > 0$ for every $\alpha$, which immediately yields $\mu(\alpha,\beta,x) = \mu_b(\alpha,\beta,x)$. □

We are now ready to prove:

PROPOSITION 2.3.

$$\sup_{t>0} t \sup\left\{Q(\mu) : \mu \in \mathcal{P}, \sum_x x\mu(x) \leq 1/t\right\} = \tilde{b}(\lambda, h). \tag{2.14}$$

PROOF. Lemma 2.2 yields

$$Q(\mu_b^{\lambda,h}) = \sup\left\{Q(\mu) : \mu \in \mathcal{P}, \sum_x x\mu(x) = \sum_x x\mu_b^{\lambda,h}(x)\right\}, \tag{2.15}$$

and that the supremum is uniquely attained. Therefore, by recalling (2.11), we obtain

$$\begin{aligned}
&\sup_{t>0} t \sup\left\{Q(\mu) : \mu \in \mathcal{P}, \sum_x x\mu(x) \leq \frac{1}{t}\right\} \\
&= \sup_{b>0} \frac{Q(\mu_b^{\lambda,h})}{\sum_x x\mu_b^{\lambda,h}(x)} = \sup_{b>0}\left(b + \frac{\log Z(b,\lambda,h)}{\sum_x x\mu_b^{\lambda,h}(x)}\right).
\end{aligned}$$
(2.16)

Let us observe that, by Lemma 2.1, the derivative with respect to $b$ of the argument of the supremum in the right-hand side is equal to

$$(2.17) \qquad -\log Z(b)\frac{f'(b)}{f^2(b)},$$

where, as before, $f(b) = \sum_x x\mu_b(x)$. Lemma 2.2 guarantees that $f'(b) < 0$ for every $b > 0$ and therefore a point $b > 0$ such that $Z(b) = 1$ is the maximum. Such a point may not exist; this is the case if $Z(b) < 1$ for every $b > 0$, that is, $Z(0) \le 1$, and in this case the supremum is achieved in the limit $b \searrow 0$ and it takes of course the value zero. □

**3. From large deviations to the variational problem.** In this section we present a proof of Proposition 1.4. This, coupled with Proposition 2.3, yields Theorem 1.2. It will be preceded by some straightforward manipulations of the free energy and by a large deviations principle for suitable Markov processes.

3.1. *Reduction to random walk excursions.* We start by recalling that the reduced free energy $\phi_\omega(\lambda, h)$ is defined as limit as $N \to \infty$ of

$$(3.1) \quad \phi_{N,\omega}(\lambda, h) = \frac{1}{N}\log \mathbb{E}\left(\exp\left(\lambda\left(\sum_{x=1}^N (\omega_x + h)\operatorname{sign}(S_x)\right) - \lambda h N\right)\right).$$

For the arguments in this section there is some technical advantage in using instead

$$(3.2) \quad \tilde{\phi}_{N,\omega}(\lambda, h) = \frac{1}{N}\log \mathbb{E}\left(\exp\left(\lambda\left(\sum_{x=1}^N (\omega_x + h)(\operatorname{sign}(S_x) - 1)\right)\right)\right).$$

Since $|\sum_{x=1}^N \omega_x| \le T_\omega$, it is immediate to see that $|\tilde{\phi}_{N,\omega}(\lambda, h) - \phi_{N,\omega}(\lambda, h)| \le \lambda T_\omega/N$ and the two quantities are therefore equivalent.

Starting from the setup of Section 1.6, recall in particular the sequence of stopping times defined in (1.17); we introduce also $\Delta \eta_k = \eta_{k+1} - \eta_k$, $\alpha_k = [\eta_k/2]_\mathbb{S}$ and $\beta_k = \alpha_{k+1}$, for every $k = 0, 1, \ldots$, and $\ell_N = \max\{j \in \mathbb{N} \cup 0 : \eta_j \le N\}$. By exploiting the up–down symmetry of the excursions of $S$ we easily arrive at

$$(3.3) \qquad \tilde{\phi}_{N,\omega}(\lambda, h) = \frac{1}{N}\log \mathbb{E}\left[\exp\left(\sum_{i=0}^{\ell_N - 1} \varphi(\lambda \xi_{\alpha_i, \beta_i} + \lambda h \Delta \eta_i)\right) R_N\right],$$

where $\varphi(t) = \log((1 + \exp(-2t))/2)$, $t \in \mathbb{R}$, and $R_N = \exp(\varphi(\lambda \sum_{x=\eta_{\ell_N}+1}^N (\omega_x + h)))$. Since the argument of $\varphi(\cdot)$ in (3.3) is bounded below by $-\lambda \xi_\star$, we may redefine $\varphi(\cdot)$ by $\varphi(\cdot) \vee \varphi(-\lambda \xi_\star)$ and therefore now $\|\varphi\|_\infty < \infty$.

3.2. *A Donsker–Varadhan large deviations principle.* Recall from Section 1.6 the definition of $\mathcal{P}$, which is a subset of the probability measures on $\mathbb{S} \times \mathbb{S} \times 2\mathbb{N}$: the latter space is endowed with the discrete topology and $\mathcal{P}$ is endowed with the topology of weak convergence.

For $m \in \mathbb{N}$ introduce the empirical measure

$$\mathbf{L}_m(\alpha, \beta, x) = \frac{1}{m} \sum_{j=0}^{m-1} \mathbf{1}_{\{\alpha_j, \beta_j, \Delta\eta_j\}}(\alpha, \beta, x). \tag{3.4}$$

PROPOSITION 3.1. *A full large deviations principle with rate functional $I$, defined as in (1.20), holds for the sequence of empirical measures $\{\mathbf{L}_m\}_m$. More explicitly, we have*

$$\begin{aligned}
-\inf_{\mu \in A^\circ} I(\mu) &\le \liminf_{m \to \infty} \frac{1}{m} \log \mathbb{P}(\mathbf{L}_m \in A) \\
&\le \limsup_{m \to \infty} \frac{1}{m} \log \mathbb{P}(\mathbf{L}_m \in A) \le -\inf_{\mu \in \overline{A}} I(\mu),
\end{aligned} \tag{3.5}$$

*where $A^\circ$ and $\overline{A}$ are, respectively, the interior and the closure of the set $A$.*

PROOF. This result is implicitly contained in [2]: on page 97 one finds the definition of the functional $\widetilde{I}$ of which $I$, that we have introduced in (1.20), is a particular case. Notice, moreover, that ([2], Lemma 2.1) provides the link between $\widetilde{I}$ and the standard (variational) Donsker–Varadhan expression.

Nevertheless, we sketch here a proof: let us observe that $(\alpha_k, \Delta\eta_{k-1})_{k=1,2,\ldots}$ is a uniformly ergodic Markov chain on $\mathbb{S} \times 2\mathbb{N}$ with transition probabilities $q((\alpha, x), (\alpha', x')) \stackrel{\text{def}}{=} p_{\alpha, \alpha'} K_{\alpha, \alpha'}(x')$. So (see [5], Chapter 6), the empirical measure $(1/m) \sum_{k=1}^{m} \delta_{(\alpha_k, \Delta\eta_{k-1}, \alpha_{k+1}, \Delta\eta_k)}$ satisfies a strong large deviations principle with rate function

$$\hat{I}(\nu) \stackrel{\text{def}}{=} \begin{cases} H(\nu | \nu_1 \otimes q), & \text{if } \nu_1 = \nu_2, \\ \infty, & \text{otherwise}, \end{cases} \tag{3.6}$$

$\nu$ a probability measure on $(\mathbb{S} \times 2\mathbb{N})^2$, $\nu_1, \nu_2$ being again the two marginals on $\mathbb{S} \times 2\mathbb{N}$ and $H$ the usual relative entropy. With standard notation we have used $\nu_1 \otimes q(\alpha, x, \alpha', x') = \nu_1(\alpha, x) q((\alpha, x), (\alpha', x'))$. Write $\hat{\nu}$ for the marginal on $\mathbb{S} \times \mathbb{S} \times 2\mathbb{N}$ given by $\hat{\nu}(\alpha, \beta, x) \stackrel{\text{def}}{=} \sum_z \nu(\alpha, z, \beta, x)$. By the contraction principle, $\mathbf{L}_m$ satisfies a large deviations principle with rate function $I_0(\mu) \stackrel{\text{def}}{=} \inf\{\hat{I}(\nu) : \hat{\nu} = \mu\}$, and it remains to show that $I_0 = I$. Evidently, $I_0(\mu) = \infty$ for $\mu \notin \mathcal{P}$. For a given $\mu \in \mathcal{P}$, set

$$\bar{\mu}(\alpha, x, \alpha', x') \stackrel{\text{def}}{=} \frac{\sum_\gamma \mu(\gamma, \alpha, x)}{\sum_{\gamma, x} \mu(\gamma, \alpha, x)} \mu(\alpha, \alpha', x'), \tag{3.7}$$

which evidently satisfies $\bar\mu_1 = \bar\mu_2$ if $\mu \in \mathcal{P}$. Now, for any $\nu$ satisfying $\nu_1 = \nu_2$ and $\hat\nu = \mu$ one has by an elementary computation $\hat I(\nu) - \hat I(\bar\mu) = H(\nu|\bar\mu) \geq 0$, and therefore $I_0(\mu) = \hat I(\bar\mu) = I(\mu)$.  □

PROOF OF PROPOSITION 1.4: UPPER BOUND. Since $1/2 < R_N \leq (1 + \exp(2\lambda \xi_\star))/2$, we may safely get rid of $R_N$ in the expression in (3.3) and it suffices to prove the statement for

$$(3.8) \qquad \hat\phi_{N,\omega}(\lambda, h) = \frac{1}{N} \log \mathbb{E}\left[\exp\left(\sum_{i=0}^{\ell_N - 1} \varphi(\lambda \xi_{\alpha_i,\beta_i} + \lambda h \Delta \eta_i)\right)\right],$$

rather than for $\tilde\phi_{N,\omega}(\lambda, h)$.

Let us first control $\limsup_N \hat\phi_{N,\omega}(\lambda, h)$. We proceed with a discretization procedure: choose a large integer $K$ and assume for simplicity that $2K$ divides $N$: we have

$$\begin{aligned}
\hat\phi_{N,\omega}(\lambda, h) &= \frac{1}{N} \log \sum_{j=0}^{K-1} \mathbb{E}\left[\exp\left(\sum_{i=0}^{\ell_N - 1} \varphi(\lambda \xi_{\alpha_i,\beta_i} + \lambda h \Delta \eta_i)\right); \right. \\
& \hspace{6cm} \left. \ell_N \in \left(\frac{jN}{2K}, \frac{(j+1)N}{2K}\right]\right], \\
(3.9) \quad &\leq \frac{1}{N} \log \sum_{j=0}^{K-1} \mathbb{E}\left[\exp\left(\sum_{i=0}^{jN/2K} \varphi(\lambda \xi_{\alpha_i,\beta_i} + \lambda h \Delta \eta_i)\right) \exp\left(\frac{N\|\varphi\|_\infty}{2K}\right); \right. \\
& \hspace{6cm} \left. \ell_N \in \left(\frac{jN}{2K}, \frac{(j+1)N}{2K}\right]\right],
\end{aligned}$$

from which we obtain that

$$(3.10) \quad \begin{aligned}
&\limsup_{N \to \infty} \hat\phi_{N,\omega}(\lambda, h) \\
&\leq \frac{\|\varphi\|_\infty}{2K} \\
&\quad \vee \max_{j=1,\dots,K-1} \limsup_{N \to \infty} \frac{1}{N} \log \mathbb{E}\left[\exp\left(\sum_{i=0}^{jN/2K} \varphi(\lambda \xi_{\alpha_i,\beta_i} + \lambda h \Delta \eta_i)\right); \right. \\
& \hspace{9cm} \left. \sum_{i=0}^{jN/2K} \Delta \eta_i \leq N\right].
\end{aligned}$$

But the term in the right-hand side can be easily expressed as a functional of the empirical measure, namely:

$$q := \limsup_{N \to \infty} \frac{1}{N} \log \mathbb{E}\left[\exp\left(\sum_{i=0}^{jN/2K} \varphi(\lambda \xi_{\alpha_i,\beta_i} + \lambda h \Delta \eta_i)\right); \sum_{i=0}^{jN/2K} \Delta \eta_i \leq N\right]$$

$$
\begin{aligned}
&= \limsup_{N \to \infty} \frac{1}{N} \log \mathbb{E}\Bigg[\exp\bigg(N\bigg(\frac{j}{2K}\bigg) \\
&\qquad\qquad\qquad\qquad \times \sum_{\alpha,\beta,x} \mathbf{L}_{jN/2K}(\alpha,\beta,x)\varphi(\lambda\xi_{\alpha,\beta} + \lambda h x)\bigg); \\
&\qquad\qquad\qquad\qquad\qquad \sum_{\alpha,\beta,x} x\mathbf{L}_{jN/2K}(\alpha,\beta,x) \leq \frac{2K}{j}\Bigg].
\end{aligned}
$$
(3.11)

Since $\{\mu : \sum_{\alpha,\beta,x} x\mu(\alpha,\beta,x) \leq c\}$ is a closed set, we may apply the upper bound of the large deviations principle in Proposition 3.1 to obtain

$$
q \leq \sup\Bigg\{\frac{j}{2K}\bigg(\sum_{\alpha,\beta,x} \mu(\alpha,\beta,x)\varphi(\lambda\xi_{\alpha,\beta} + \lambda h x) - I(\mu)\bigg) : \\
\sum_{\alpha,\beta,x} x\mu(\alpha,\beta,x) \leq \frac{2K}{j}\Bigg\}.
$$
(3.12)

Recalling (3.10) and taking $K \to \infty$, we arrive at

$$
\limsup_{N \to \infty} \hat{\phi}_{N,\omega}(\lambda, h)
$$
(3.13)
$$
\leq \sup_{t \in (0,1/2]} \sup_{\substack{\mu \in \mathcal{P} : \\ \sum_x x\mu(x) \leq 1/t}} t\bigg(\sum_{\alpha,\beta,x} \mu(\alpha,\beta,x)\varphi(\lambda\xi_{\alpha,\beta} + \lambda h x) - I(\mu)\bigg),
$$

which proves that the left-hand side in (1.22) is not larger than the right-hand side. □

PROOF OF PROPOSITION 1.4: LOWER BOUND. As in the previous proof we may concentrate on $\hat{\phi}_{N,\omega}(\lambda, h)$, see (3.8), and on the process $\{\Delta\eta_j\}_j$ rather than on the whole $S$-path. For $b > 0$ we consider the measure $\mu_b = \mu_b^{\lambda,h}$, defined in (2.2), which naturally defines the Markov process $\{\alpha_i, \beta_i, \Delta\eta_i\}_i$ introduced in the proof of Lemma 2.1. For definiteness we consider the stationary process conditioned to $\alpha_0 = 0$ and denote by $\mathbb{P}_b^{(N)}$ the law of $\{\ell_N, \{\Delta\eta_j\}_{j=0,\ldots,\ell_N-1}\}$. Notice that if $b = 0$ and $\lambda = 0$, then the process $\{\alpha_i, \beta_i, \Delta\eta_i\}_i$ is the one associated to the simple random walk $S$. We denote by $\mathbb{P}^{(N)}$ the measure $\mathbb{P}_0^{(N)}$ with $\lambda = 0$.

By applying the Jensen inequality we have

$$
(3.14) \quad \hat{\phi}_{N,\omega}(\lambda, h) \geq \mathbb{E}_b^{(N)}\Bigg[\frac{1}{N}\sum_{i=0}^{\ell_N-1} \varphi(\lambda\xi_{\alpha_i,\beta_i} + \lambda h \Delta\eta_i)\Bigg] - \frac{1}{N}H(\mathbb{P}_b^{(N)}|\mathbb{P}^{(N)}),
$$

in which $H$ still denotes the relative entropy. The relative entropy term can be evaluated directly and one obtains

$$
\begin{aligned}
&\liminf_{N\to\infty} \hat{\phi}_{N,\omega}(\lambda, h) \\
&\qquad \geq \liminf_{N\to\infty} \mathbb{E}_b\left[\frac{1}{N} \sum_{i=0}^{\ell_N-1} \varphi(\lambda \xi_{\alpha_i,\beta_i} + \lambda h \Delta \eta_i)\right] \\
&\qquad \quad - \limsup_{N\to\infty} \mathbb{E}_b\left[\frac{1}{N} \sum_{i=0}^{\ell_N-1} G(\alpha_i, \beta_i, \Delta \eta_i)\right],
\end{aligned}
\tag{3.15}
$$

where $G(\alpha, \beta, x) \equiv \log(\mu_b(\alpha, \beta, x)/\mu_b(\alpha)p_{\alpha,\beta}K_{\alpha,\beta}(x)) = \Phi_{\alpha,\beta}^{\lambda,h}(x) - bx - \log Z + \log(v_\alpha/v_\beta)$. By ergodicity $\lim_N \ell_N/N = 1/\sum_x x\mu_b(x)$ $\mathbb{P}_b$-a.s.; this implies that $\sum_{i=0}^{\ell_N-1} G(\alpha_i, \beta_i, \Delta \eta_i)/N$ converges a.s. to $I(\mu_b)/\sum_x x\mu_b(x)$. Since $G(\alpha, \beta, x)$ is bounded above, by Fatou's lemma one obtains that the superior limit of $\mathbb{E}_b[\sum_{i=0}^{\ell_N-1} G(\alpha_i, \beta_i, \Delta \eta_i)/N]$ is bounded above by $I(\mu_b)/\sum_x x\mu_b(x)$. By recalling that $\varphi(\cdot)$ is bounded and by applying once again the ergodic theorem one obtains

$$
\begin{aligned}
\liminf_{N\to\infty} \hat{\phi}_{N,\omega}(\lambda, h) &\geq \frac{\sum_{\alpha,\beta,x} \varphi(\lambda \xi_{\alpha,\beta} + \lambda h x) \mu_b(\alpha, \beta, x) - I(\mu_b)}{\sum_x x\mu_b(x)} \\
&= \frac{Q(\mu_b)}{\sum_x x\mu_b(x)}.
\end{aligned}
\tag{3.16}
$$

Set $t = 1/\sum_x x\mu_b(x)$; Lemma 2.2 allows then to replace in the previous expression $\mu_b$ with any $\mu \in \mathcal{P}$ such that $\sum_x x\mu(x) = 1/t$. Optimizing over $t$ leads to the desired lower bound and proof of Theorem 1.4 is complete. $\square$

**4. Existence and monotonicity of the critical line.** We now go through some soft arguments that yield a partial proof of Proposition 1.1. We prove the following:

LEMMA 4.1. *There exists a nondecreasing function $h_c(\cdot):[0,\infty) \to [0,1]$, continuous if its domain is restricted to $(0,\infty)$, such that $\mathcal{D} = \{(\lambda, h): h \geq h_c(\lambda)\}$.*

REMARK 4.2. In order to fill the gap with Proposition 1.1 one needs to show that $h_c(\lambda)$ vanishes as $\lambda \searrow 0$ and that it tends to 1 as $\lambda \nearrow 0$, as well as the fact that the image is $[0, 1)$ rather than $[0, 1]$. To establish this one needs some quantitative bounds on $\phi_\omega(\cdot)$. Theorem 1.3 of course largely provides the needed bound for small $\lambda$; for large $\lambda$ we refer to Section 5.2.

PROOF OF LEMMA 4.1. We first collect some elementary facts:

1. $\phi_\omega(0,h) = 0$ for every $h$ and $\phi_\omega(\lambda, h) \geq 0$ for every $\lambda$ and $h$.
2. $\phi_\omega(\cdot, h)$ is a convex function for every $h$. This simply follows from the convexity in $\lambda$ of $\log Z_{N,\omega}$, see (1.4).
3. $\phi_\omega(\lambda, \cdot)$ is nonincreasing, besides being convex (proven as in point 2) and therefore continuous, for every $\lambda$. This follows from (1.22) and the fact that $\widetilde{\Phi}_{\alpha,\beta}^{\lambda,\cdot}(x)$ is nonincreasing for every $\alpha, \beta, \lambda$ and $x$.
4. If $(\lambda, h) \in \mathcal{D}$ and if $\Delta\lambda$ and $\Delta h$ are two positive numbers such that $\Delta h \geq \Delta\lambda(1-h)/\lambda$, then $(\lambda + \Delta\lambda, h + \Delta h) \in \mathcal{D}$. This follows by (1.22) once we observe that $\widetilde{\Phi}_{\alpha,\beta}^{\lambda,h}(x) = \log(1 + \exp(-\lambda(\xi_{\alpha,\beta} + hx))) - \log 2$ is a decreasing function of $\lambda(\xi_{\alpha,\beta} + hx)$. By using that $\xi_{\alpha,\beta} \geq -x$ one directly verifies that choosing $\Delta\lambda$ and $\Delta h$ as above implies $(\lambda + \Delta\lambda)(\xi_{\alpha,\beta} + (h + \Delta h)x) \geq \lambda(\xi_{\alpha,\beta} + hx)$.
5. If $h \geq 1$, then $\phi_\omega(\cdot, h)$ is nonincreasing. This again follows from the fact that $\widetilde{\Phi}_{\alpha,\beta}^{\lambda,h}(x)$ is a decreasing function of $\lambda(\xi_{\alpha,\beta} + hx)$: if $h \geq 1$, then $\xi_{\alpha,\beta} + hx \geq 0$, so that $\widetilde{\Phi}_{\alpha,\beta}^{\cdot,h}(x)$ is nonincreasing and (1.22) implies the result.

Let us use these five facts (*points*): first of all we observe that point 3 guarantees the existence of $h_c(\cdot)$ and that $h_c(\lambda) = \inf\{h : \phi_\omega(\lambda, h) = 0\}$. By points 1 and 2 we have that $\phi_\omega(\cdot, h)$ is nondecreasing, so $h_c(\cdot)$ is nondecreasing. Points 1 and 5 guarantee that $h_c(\lambda) \leq 1$ for every $\lambda \geq 0$. Finally, point 4 implies the continuity of $h_c(\cdot)$ on the positive semi-axis because it implies that the incremental ratios of $h_c(\cdot)$ at the point $\lambda$ are bounded above by $(1 - h_c(\lambda))/\lambda$. $\square$

## 5. Asymptotics for small and large $\lambda$: the proof of Theorem 1.3.

5.1. *Small $\lambda$ asymptotics.* Set $h = m\lambda^3$, $m$ a positive number, and write $\Phi_{\alpha,\beta}(\cdot) = \Phi_{\alpha,\beta}^{\lambda, m\lambda^3}(\cdot)$: we choose to work with the latter defined as in (1.23). Set also $\widetilde{A}(\lambda) = A(0, \lambda, m\lambda^3)$. We write

$$(5.1) \qquad \widetilde{A}_{\alpha,\beta}(\lambda) = p_{\alpha,\beta} + p_{\alpha,\beta} \sum_x K_{\alpha,\beta}(x)[\exp(\Phi_{\alpha,\beta}(x)) - 1],$$

and we decompose the second term in the right-hand side in three terms:

$$(5.2) \qquad \begin{aligned} \sum_x & K_{\alpha,\beta}(x)[\exp(\Phi_{\alpha,\beta}(x)) - 1] \\ &= \sum_{x \leq \lambda^{-5/2}} K_{\alpha,\beta}(x)[\exp(\Phi_{\alpha,\beta}(x) + m\lambda^4 x) - 1]\exp(-m\lambda^4 x) \\ &\quad + \sum_{x \leq \lambda^{-5/2}} K_{\alpha,\beta}(x)[\exp(-m\lambda^4 x) - 1] \\ &\quad + \sum_{x > \lambda^{-5/2}} K_{\alpha,\beta}(x)[\exp(\Phi_{\alpha,\beta}(x)) - 1] \end{aligned}$$

$$:= T_1 + T_2 + T_3.$$

By using the fact that for $x \ll \lambda^{-3}$

$$\text{(5.3)} \qquad \frac{\exp(\Phi_{\alpha,\beta}(x) + m\lambda^4 x) - 1}{\lambda^2} = \frac{\xi_{\alpha,\beta}^2}{2}(1 + o(1)),$$

one easily sees that $\lim_{\lambda \searrow 0} \lambda^{-2} T_1 = \xi_{\alpha,\beta}^2/2$. Moreover, by using that $1 - e^{-t} \leq t$ for $t \geq 0$, one directly obtains that $T_2 = O(\lambda^{11/4}) = o(\lambda^2)$. For $T_3$ we observe that

$$T_3 = \sum_{x > \lambda^{-5/2}} K_{\alpha,\beta}(x) \exp(\Delta)[\exp(\psi(m\lambda^4 x) - m\lambda^4 x) - 1]$$

$$\text{(5.4)} \qquad + \sum_{x > \lambda^{-5/2}} K_{\alpha,\beta}(x)[\exp(\Delta) - 1]$$

$$= T_4 + T_5,$$

with $\Delta = \psi(\lambda \xi_{\alpha,\beta} + m\lambda^4 x) - \psi(m\lambda^4 x)$. Taylor's expansion yields immediately that $|\Delta| \leq \xi_\star \lambda$. This suffices to show that $T_5$ is inessential:

$$\text{(5.5)} \qquad |T_5| \leq K \sum_{x > \lambda^{-5/2}} x^{-3/2} |e^\Delta - 1| \leq c\lambda^{1+1/4} \lambda = o(\lambda^2).$$

For $T_4$ we approximate the sum by an integral using (1.10) and obtain thus

$$\lim_{\lambda \searrow 0} \lambda^{-2} \sum_{x > \lambda^{-5/2}} K_{\alpha,\beta}(x)[\cosh(m\lambda^4 x) \exp(-m\lambda^4 x) - 1]$$

$$\text{(5.6)} \qquad = \frac{\mathsf{C}_K \sqrt{m}}{2 T_\omega p_{\alpha,\beta}} \int_0^\infty \frac{1}{r^{3/2}} [\cosh(r) \exp(-r) - 1] \, dr$$

$$= -\frac{\mathsf{C}_K \sqrt{2\pi m}}{2 T_\omega p_{\alpha,\beta}}.$$

The conclusion is that

$$\text{(5.7)} \qquad \widetilde{A}_{\alpha,\beta}(\lambda) = p_{\alpha,\beta} + \lambda^2 \left( p_{\alpha,\beta} \frac{\xi_{\alpha,\beta}^2}{2} - \frac{\mathsf{C}_K \sqrt{2\pi m}}{2 T_\omega} \right) + o(\lambda^2).$$

By applying, for example, (2.1) one directly writes the corresponding expansion for the maximal eigenvalue $\widetilde{Z}(\lambda)$ of $\widetilde{A}(\lambda)$ (recall that $p_{\alpha,\beta}$ is bi-stochastic):

$$\text{(5.8)} \qquad \widetilde{Z}(\lambda) = 1 + \lambda^2 \left( \frac{1}{T_\omega} \sum_{\alpha,\beta} p_{\alpha,\beta} \frac{\xi_{\alpha,\beta}^2}{2} - \mathsf{C}_K \sqrt{\frac{\pi}{2} m} \right) + o(\lambda^2).$$

In view of Theorem 1.2, formula (1.14) is proven with $m_\omega$ equal to the value of $m$ for which the term between brackets in (5.8) is zero.

5.2. *Large $\lambda$ asymptotics.* We set $h = 1 - (M/\lambda)$, $M > 0$, and the argument we are going to use is based on the observation that if we define for $[x]_{\mathbb{S}} = \beta - \alpha$

$$(5.9) \quad \widehat{\Phi}_{\alpha,\beta}(x) = -\log 2 + \begin{cases} 0, & \text{if } \xi_{\alpha,\beta} > -x, \\ \log(1 + \exp(2Mx)), & \text{if } \xi_{\alpha,\beta} = -x, \end{cases}$$

then

$$(5.10) \qquad 0 \le \Phi_{\alpha,\beta}^{\lambda, 1-(M/\lambda)}(x) - \widehat{\Phi}_{\alpha,\beta}(x) \le \exp(-4\lambda),$$

with $\Phi$ as in (1.11).

We have therefore that $Z(0, \lambda, 1 - (M/\lambda))$ tends to $\widehat{Z}(M)$, the principal eigenvalue of $p_{\alpha,\beta} \sum_x K_{\alpha,\beta}(x) \exp(\widehat{\Phi}_{\alpha,\beta}(x))$. Notice that $\widehat{Z}(\cdot)$ is increasing, that $\widehat{Z}(0) = 1/2$ and that $\widehat{Z}(M)$ tends to infinity as $M \nearrow \infty$ (this can be seen, e.g., by applying Theorem 1.4 in [10], Chapter 2) so that $M_\omega = \widehat{Z}^{-1}(1)$. This proves (1.16) and the proof of Theorem 1.3 is complete.

**Acknowledgment.** We are grateful to Massimiliano Gubinelli for having pointed out the result in Lemma 2.1.

INSTITUT FÜR MATHEMATIK
MATHEMATISCH–NATURWISSENSCHAFTLICHE FAKULTÄT
UNIVERSITÄT ZÜRICH
WINTERTHURERSTRASSE 190, CH-8057
ZÜRICH
E-MAIL: eb@amath.unizh.ch

U.F.R. MATHEMATIQUES
UNIVERSITÉ PARIS 7–DENIS DIDEROT
CASE 7012
2 PLACE JUSSIEU 75251
PARIS CEDEX 05
FRANCE
E-MAIL: giacomin@math.jussieu.fr